\documentclass[10pt]{article}

\usepackage{graphics}
\usepackage{graphicx}

\usepackage[a4paper, left=35mm,right=35mm,top=34mm,bottom=34mm]{geometry}
\usepackage[utf8]{inputenc}
\usepackage[T1]{fontenc}
\usepackage[english]{babel}

\usepackage{enumerate}
\usepackage{graphicx}
\usepackage{hyperref}
\hypersetup{
    colorlinks=true,
    linkcolor=blue,
    filecolor=magenta,      
    urlcolor=cyan,
}

\usepackage{mathtools,amsthm,amssymb,amsfonts}
\usepackage{algorithm}
\usepackage{algorithmic}
\makeatother
\theoremstyle{plain}

\theoremstyle{definition}

\theoremstyle{remark}

\usepackage{caption} 
\usepackage[tight,footnotesize]{subfigure}

\usepackage{fancyhdr}
\usepackage{array}
\usepackage{xcolor}
\usepackage{ecp/ecpmashpc}

\author{
  {\normalsize Abal-Kassim Cheik Ahamed}\thanks{\'Ecole Centrale Paris, France.}
	\and
  {\normalsize Alban Desmaison}\thanks{\'Ecole Centrale Paris, France.}
	\and
  {\normalsize Fr\'ed\'eric Magoul\`es}\thanks{\'Ecole Centrale Paris, France
    (correspondence, frederic.magoules@hotmail.com).}
		}	
\title{Fast and Green Computing with Graphics Processing Units for solving Sparse Linear Systems}
\date{}

\begin{document}
\maketitle
\thispagestyle{fancy}

\begin{abstract}
\noindent In this paper, we aim to introduce a new perspective when comparing highly parallelized algorithms on GPU: the energy consumption of the GPU. We give an analysis of the performance of linear algebra operations, including addition of vectors, element-wise product, dot product and sparse matrix-vector product, in order to validate our experimental protocol.
We also analyze their uses within conjugate gradient method for solving the gravity equations on Graphics Processing Unit (GPU). Cusp library is considered and compared to our own implementation with a set of real matrices arrising from the Chicxulub crater and obtained by the finite element discretization of the gravity equations. The experiments demonstrate the performance and robustness of our implementation in terms of energy efficiency.
\end{abstract}


\section{Introduction}
\label{sec:introduction}

Central Processing Unit (CPU) and Graphics Processing Unit (GPU) are very different processor, one is designed to be able to do as many different task as possible with the best performance possible. The other one is specialized in intensive and repetitive work. Even though the CPUs are efficient for calculation, logic algorithm and data management, they lack the possibility to do multiple calculus at the same time and process a large amount of data quickly. This kind of task, corresponding at first to graphic processing, pushed the evolution of the GPU. They are processors made to be able to handle large amount of data with simple task. This completely different architecture allows the GPU to be way more efficient for some tasks than a CPU. Most importantly, for scientific computing, these kind of processor are a large breakthrough.
One interesting point about GPUs is that their large number of parallel computing cores allows them to run at a much lower frequency than a CPU. This specialty allows them to be less energy consuming while delivering a way bigger computing capability. Indeed the top list of supercomputers presented at~\cite{GPU:top500:2013} increasingly include GPUs to be able to increase their computing capabilities, and in particular most importantly the top list of green~\cite{GPU:green500:2013} supercomputers.
Even though we can see that the GPUs seems to have a real capability to both improve computing possibilities and reduce energy consumption of the supercomputers, few specific studies have been done to measure the exact energy consumption of a GPU and how to limit this energy consumption when writing a specific algorithm.

The rest of the paper is organized as follows. Section~\ref{sec:cuda_prog_model} introduces an overview of GPU programming model for non familar readers. Section~\ref{sec:experiments_configuration} presents the hardware configuration involved in this paper. The expermental protocol and requirements are detailed in Section~\ref{sec:protocol_experimentation}. In Section~\ref{sec:exp_results}, we report, evaluate and analyze numerical results. Finally, concluding remarks are given in Section~\ref{sec:concluding_remarks}.

\section{GPU programming model}
\label{sec:cuda_prog_model}

\subsection{History and evolution of GPU}

This section is going to present the history of GPU and to give a brief introduction of the main challenge for GPUs in High Performance Computing, i.e., Energy Consumption. GPU appeared early in the history of computer science. Indeed, when the computers started to have more and more advanced graphics, the CPU become a limit to the evolution of graphics. 
Graphics Processors were added to computers to handle the rendering and the very specific tasks associated with graphic rendering.
These processors quickly evolved to have more and more tasks to do before the final rendering of the screen.
At the end of the century, these Graphics Processors, became Graphics Processing Units that were doing even more tasks, from the image transformation (4x4 matrix multiplication) to lightning. We also start creating some standards for graphic processing like openGL and DirectX.
In the first part of the years 2000s, the quick evolution of the standards for rendering was too fast for the GPU makers. We started to see some low level programming capabilities having the same hardware and doing different functions. This technology started to give the possibility for a GPU to adapt its hardware~\cite{GPU:BYFWA:2009} to the processing that was needed.
The last step of this evolution happened in 2006, with the apparition of the General Purpose GPU, composed of only one type of processing unit that could be used for any task needed from the CPU. This new architecture allowed a better balance of the charge depending on the processing needed. It also introduced new possibilities for computation, these massive parallel processors gave scientists new possibilities to speedup their computations.
Here is a quick summary of the early stages of the Graphic Processors.

\begin{figure}[!ht]
\centering
\includegraphics[scale=0.30]{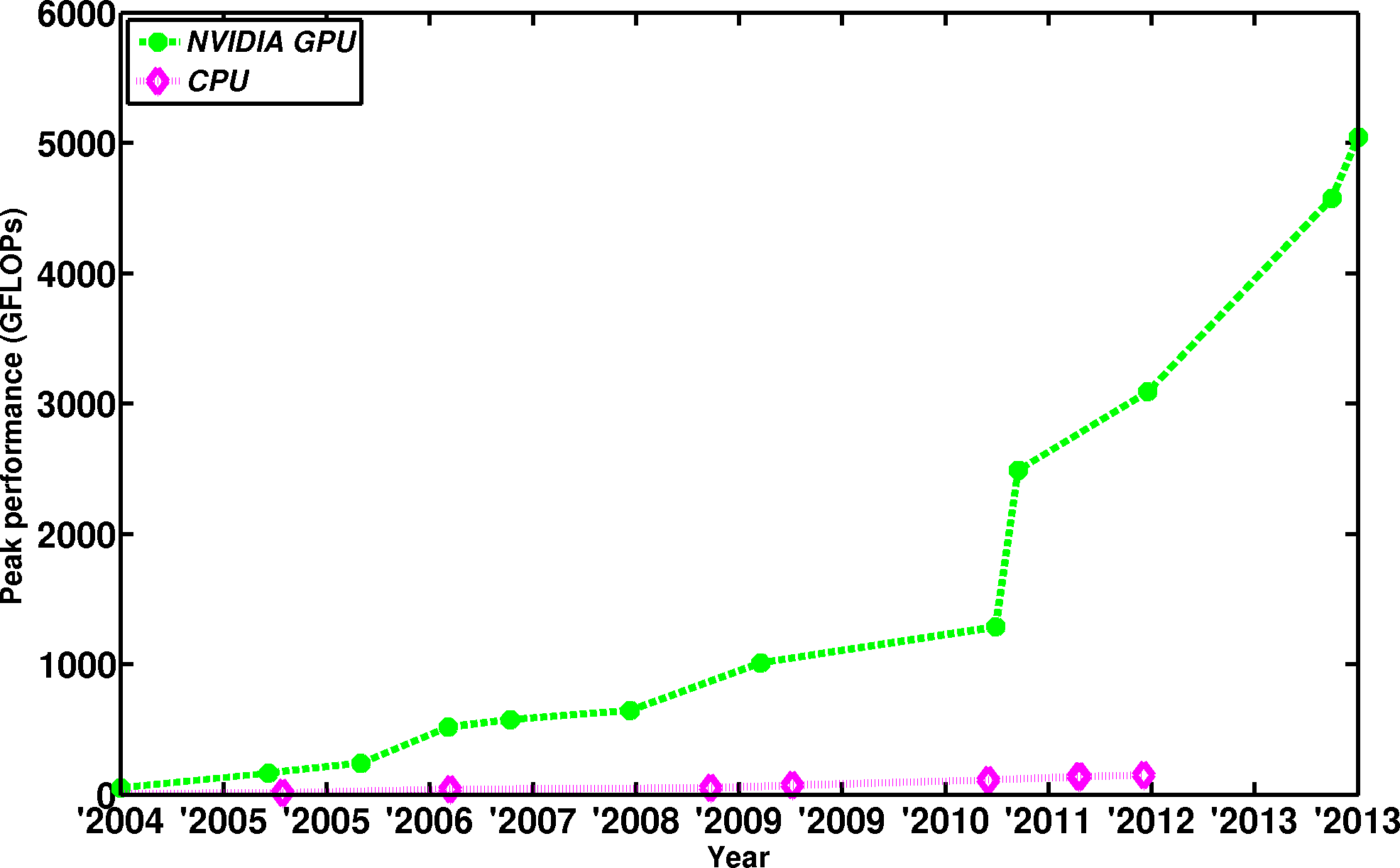}
\caption{Evolution of the peak of performance for CPU and NVIDIA GPU}
\label{fig:img:history_gpu_nvidia}
\end{figure}

Figure~\ref{fig:img:history_gpu_nvidia} reports the peak of performance of CPU and GPU from 2004 to 2013. As we can see in this figure, GPU reaches Tera floating point operations per second (Flops) when CPU hardly reaches 1GFlops. 

One of the quick advantages of the use of GPUs for scientific calculus is also its very good energetic efficiency.
When executing a properly parallelized algorithm, the GPU allows the user to drastically increase its computing capabilities. This leads to either shorter computation time or less hardware needed to run the same code with a constant time. Both of these possibilities result in a great decrease of the energy consumption to realize a calculation.
On their website, NVIDIA, explains that HESS, an oil and gas company changes a 2000 CPU cluster for a 32 Tesla S1070 (NVIDIA GPGPU specialized for High Performance Computing), keeping the same performances and reducing the energy consumption from 1.34 megaWatts to 47 kiloWatts. NVIDIA is presenting a few examples like this one on their website \url{http://www.nvidia.com}.

Even though the GPUs have been widely used for many years now, the green computing problem only appeared a few years ago, and energy efficiency on GPUs is mainly obtained by coarse grain Distributed File System (DFS) and clock gating. Currently, the GPGPUs bottleneck is located in the off-card memory bandwidth. Thus, it is not useful to have more and more powerful processing units, instead, Graphic Cards manufacturers are trying to make the current GPUs more and more efficient from an energetic point of view.

\subsection{CUDA Architecture}

One of the most widely used languages to program GPGPUs is CUDA~\cite{GPU:CUDA4.0:2011}, this language has been developed by NVIDIA to allow everyone to use their GPGPUs hardware.
CUDA is an extension of the C programming language, it is mainly a library along with a new compiler that is able to compile both CPU and GPU code depending on what the user specified.
The whole purpose of CUDA is to be able to declare a new kind of function, a Kernel that will be called from a regular CPU function. But this code will be executed on the GPU following a user specified pattern. Actually, CUDA architecture is separated in different layers:
\begin{itemize}
\item The lower level is a single scalar core that can execute one single operation, called a CUDA core.
\item A group of 32 CUDA cores is called a streaming multiprocessor, they are the basic component of the CUDA hardware architecture.
\end{itemize}

\begin{figure}[!ht]
\centering
\includegraphics[scale=0.30]{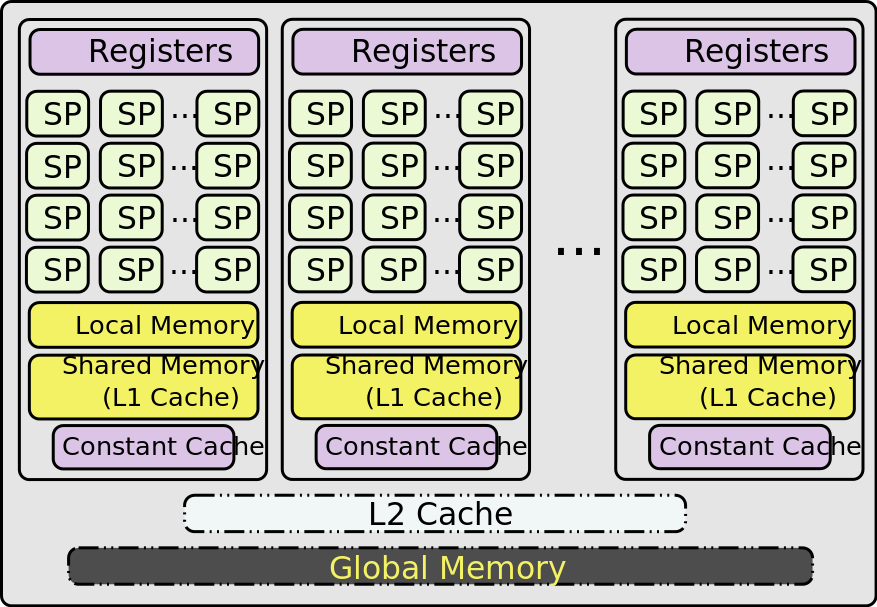}
\caption{GPU NVIDIA Fermi architecture}
\label{fig:img:gpu_fermi_archi}
\end{figure}

\begin{figure}[!ht]
\centering
\includegraphics[scale=0.30]{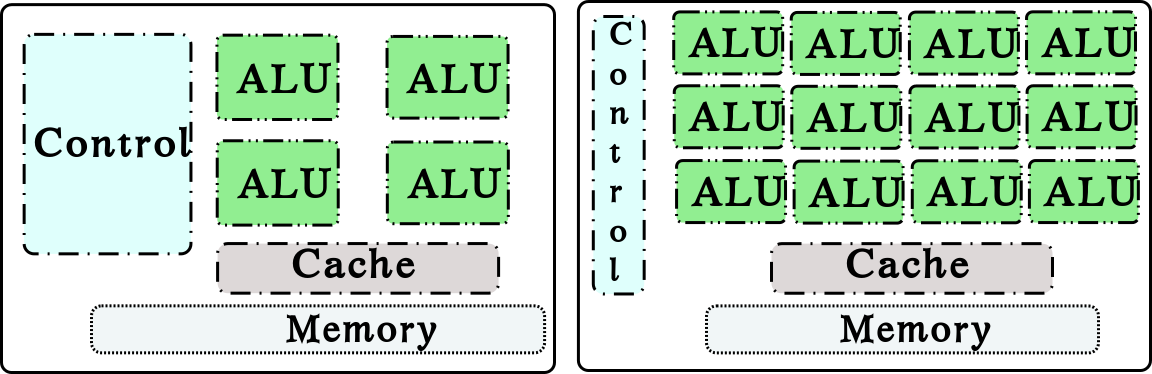}
\caption{Difference between of CPU and GPU architetcure}
\label{fig:img:cpu_gpu_architetcures}
\end{figure}

Figure~\ref{fig:img:gpu_fermi_archi} shows a simplified schema of a NVIDIA Fermi GPU. Each block represents a Streaming Multiprocessor (SM).
As we can see in Figure~\ref{fig:img:cpu_gpu_architetcures}, the main difference between CPU and GPU is due to their architecture. GPU has several ALUs.

\begin{figure}[!ht]
\centering
\includegraphics[scale=0.35]{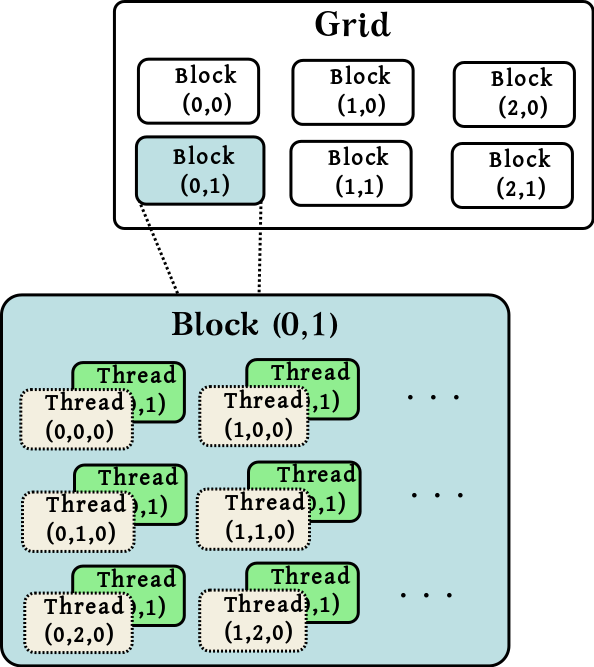}
\caption{Threading organization}
\label{fig:img:grid_block_v5}
\end{figure}

\begin{figure}[!ht]
\centering
\includegraphics[scale=0.35]{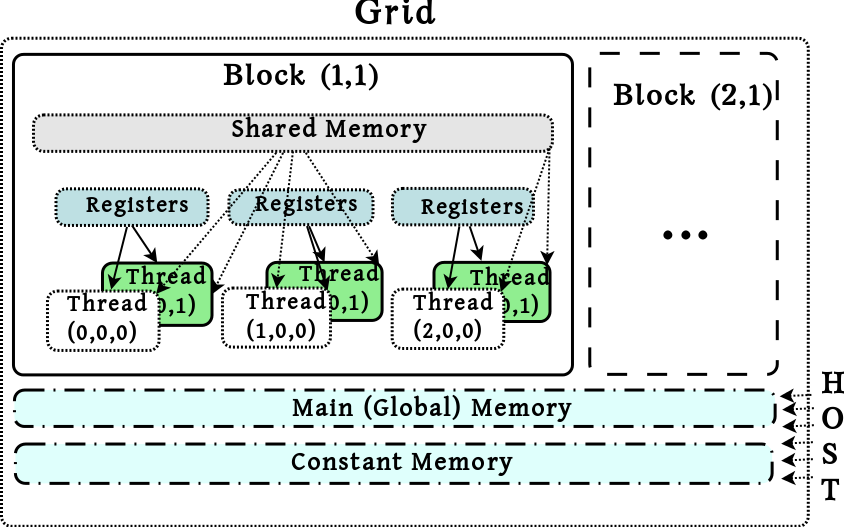}
\caption{CUDA memory level}
\label{fig:img:memory_v5}
\end{figure}

This hardware architecture is used to launch multiple instances of the same application at the same time using the Single Instruction Multiple Data (SIMD) architecture. To facilitate the management of these applications that are launched in parallel, CUDA is creating an architecture for them:
\begin{itemize}
\item Each instance of the application is called a Thread. Each thread is executed on one CUDA core.
\item Threads are gathered in Thread Blocks that can hold a certain amount of threads using a three dimensional localization of each thread: each thread is placed in a Block from its coordinates (x, y, z)
\item Thread Blocks are forming the Grid which are the biggest group of threads in CUDA. The Blocks are also located with three dimensional coordinates in the grid.
\end{itemize}
The threading organization is summarized in Figure~\ref{fig:img:grid_block_v5}.
When a program launches a kernel, it has to choose the number of Blocks in the Grid and the number of threads per blocks. Of course, these values are very important to get the best performances.
Another very important element for the performance of a CUDA kernel is the concept of warp. A warp is a set of 32 Threads. A warp is the smallest group of thread managed by the GPU at runtime. Each warp will be executed on one multiprocessor and only one instruction is executed at a time for the whole warp. It is very important to have all the threads in a warp executing the same instructions at the same time, or some of the threads will just stall, waiting for the others.
The memory architecture of cuda is very similar to a CPU one: global memory and a set of different level of cache. The different level of main memory in CUDA is given in Figure~\ref{fig:img:memory_v5}. 

\begin{itemize}
\item The constant memory is a small memory that can be used to broadcast some elements to multiple SM at the same time. It allows a bigger bandwidth in the case that multiple threads access the same data at the same time and does not modify its value. The user can write in this memory directly from the CPU.
\item The global memory is the main memory of the GPU, it is very similar to the CPU RAM. The user can write in this memory directly from the CPU.
\item The shared memory, it is between the Global Memory and the registers. This small amount of memory can be accessed by any thread inside a block. Shared memory is typically 16KB on NVIDIA GPUs and can be seen as explicitly managed cache shared inside a Streaming Multiprocessor.
\item The registers are the memory associated with each thread, it cannot be shared with other threads.
\end{itemize}

When working with more details of the architecture, you can see that there are two level of caches between the registers and the global memory. An L2 level that is shared through all SM and an L1 level inside each SM.
This memory architecture is very important because as we said before, the memory bandwidth is the bottleneck for GPU computing. They are many ways to limit the time spent to read and write in memory by using both the constant memory for constant variables and the shared memory to get some data once from the global memory and then share it through the whole block.

\subsection{Sparse matrix formats}
\label{subsec:sparse-matrix_formats}

\begin{figure}[!ht]
\centering
\includegraphics[scale=0.28]{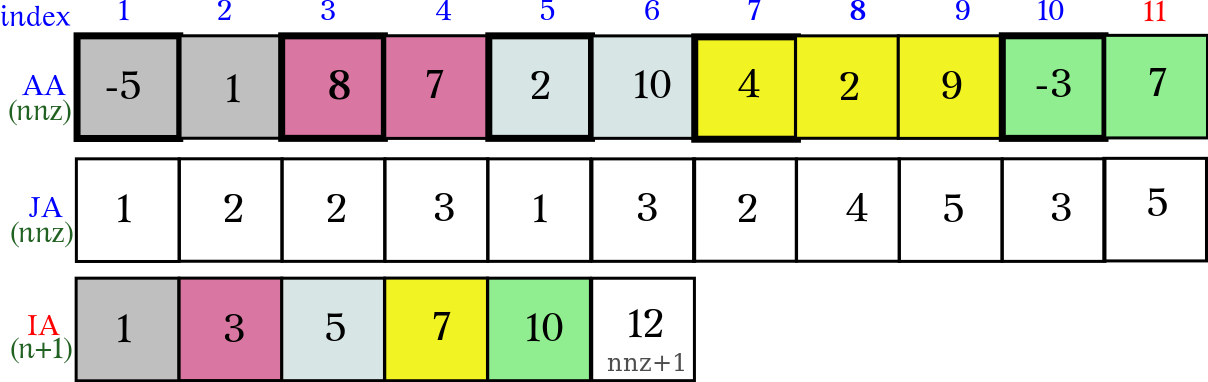}
\caption{Compressed sparse row storage format (CSR) of matrix Table~\ref{tab:matrix_and_pattern}}
\label{fig:img:csr_format}
\end{figure}

A key element in these computations are sparse matrix. These matrix are usually really large and contains mainly zero coefficients. They appears when we use numerical methods to solve partial differential equations. If their non-zero elements are creating a regular pattern along diagonals, they are called \emph{structured}, \emph{unstructured} otherwise.
The performances of the algorithms that we use can be strongly modified by the structure of the matrix \cite{cheikahamed:2012:inproceedings-2}.
Because of their size, these matrix cannot be stored as simple set of vectors. Instead, taking advantage of their great amount of zero values, several data storage structures for these matrices has been developed, for example, the Compressed-Sparse Row (CSR), Coordinate (COO), ELLPACK (ELL), Hybrid (HYB), etc.
The data storage structure that we use in this work is the CSR format. We mainly chose this one because it is the most popular and allows an easy comparison of our results with existing work.
Fig.~\ref{fig:img:csr_format} presents the storage structure used in CSR: we use three one dimensional arrays, two of size $nnz$ called $AA$ and $JA$. The first one stores all the non zero elements row-wise, the second one stores the column indices of each non-zero values: $JA(k)$ is the column number of $AA(k)$ in the initial matrix $A$. The third vector, of size $n+1$ stores the indices in the two other matrices where each row begin: in the $AA$ and $JA$ vectors, the $i-th$ row starts at $IA(i)$ and finishes at $IA(i+1) - 1$.
The main advantages of CSR when used for matrix vector multiplication is a sequential access to the data with precomputed indices that allows high efficiency on memory accesses~\cite{cheikahamed:2012:inproceedings-2,cheikahamed:2013:inproceedings-3,ref:cheikahamed:2012:inproceedings-1:DARG:2013}. The main tradeoff is that this format ignores any dense substructure that could be handled more efficiently. Moreover, when accessing the vector, we need to use $JA$ to get the indices. This kind of indirect indexing should be avoided when possible. CSR storage of matrix $A$ (Table~\ref{tab:matrix_and_pattern}) is presented in Figure~\ref{fig:img:csr_format}.
\begin{table}[!ht]
\centering
\begin{tabular}{m{7cm}m{4cm}}
\scalebox{1.4}{
$A=\begin{pmatrix}
\colorbox{gray}{-5} & \colorbox{gray}{1} & 0 & 0 & 0 \\
 0 & \colorbox{gray}{8} & \colorbox{gray}{7} & 0 & 0 \\
 \colorbox{gray}{2} & 0 & \colorbox{gray}{10} & 0 & 0 \\
 0 & \colorbox{gray}{4} & 0 & \colorbox{gray}{2} & \colorbox{gray}{9} \\
 0 & 0 & \colorbox{gray}{-3} & 0 & \colorbox{gray}{7}
\end{pmatrix}$
}
&
\includegraphics[scale=0.3]{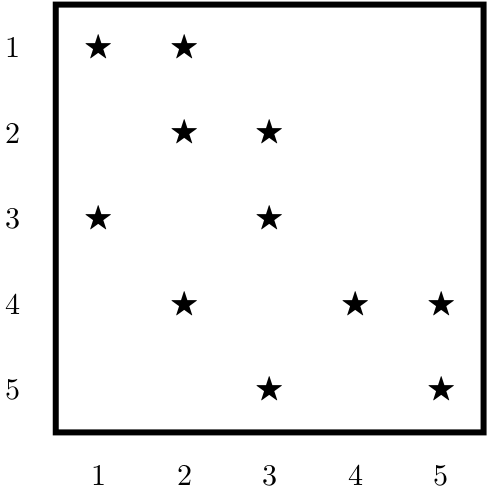}\\
\end{tabular}
\caption{Left (matrix), Right (matrix pattern)}
\label{tab:matrix_and_pattern}
\end{table}

\section{Hardware configuration}
\label{sec:experiments_configuration}

For all the experiments that we performed, we used a standard computer with an i7 920 2.67Ghz (quad core processor) and 12GB of RAM.
The Graphic cards that we used were GTX275 NVIDIA cards with 89MB of memory and Tesla K20 with 12GB memory.
Both of these GPUs have simple and double precision capability.
During the experiment, to avoid noise as much as possible, we used two GTX275. One was doing the rendering on the screen and the other one was used for the experiment.
{\bf The main purpose of our hardware setup was to be able to measure exactly the energy consumption~\cite{LDHDBX2011} of the GPU without the consumption of the other components of the computer}.
The GPU is powered from different points:
\begin{itemize}
\item The direct cables from the computer alimentation, usually 8pins or 2x8pins for the most power consuming cards. This alimentation is a 12V power-supply.
\item From the PCIe~\cite{Budruk:2003:PES:861280} itself, there are some pins that supply 12V power-supply to the GPU.
\item From the PCIe, there are some pins that supply 3.3V power-spply to the GPU.
\end{itemize}

\begin{figure}[!ht]
\centering
\includegraphics[scale=0.12]{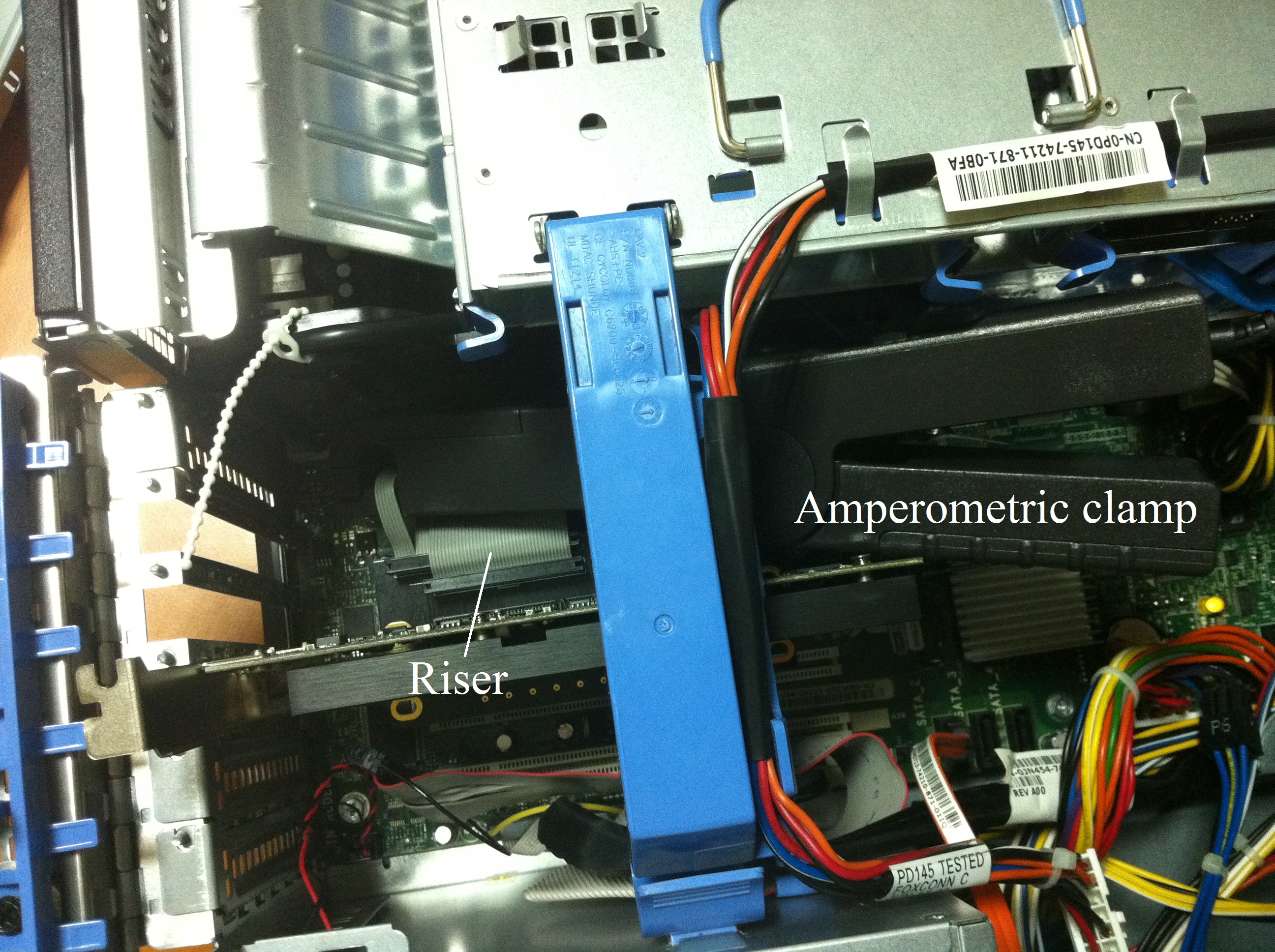}
\caption{Clamp on computer}
\label{fig:img:clamp_on_computer}
\end{figure}

These three elements are DC alimentation. To be able to know the exact energy that the GPU is using, we need to know the current that is going through all these cables and pins.
To get the current through the wires is pretty easy using amperometric clamps. These clamps output a voltage (readable on an oscilloscope) proportional to the relative sum of all the current going through the clamp.
The main problem is that in the current setup, it is impossible to measure the current going through the PCIe pins. That is why we had to use a PCIe Riser. This element is just a connector band that allows the user to move any device connected via PCIe at a different position than directly plugged on the motherboard.
We were then able to cut this band to separate some of the wires (the ones connected to the pins that were bringing power to the GPU) to be able to put an amperometric clamp on them.
Figure~\ref{fig:img:clamp_on_computer} shows how to measure the total current given to the GPU through the PCIe port for one tension.
Finally, for our setup, we needed to get the data from the amperometric clamps to the computer to be able to correlate the program that we where launching with the actual power consumption. For this purpose, we used an Ethernet-connected oscilloscope that is able to share its measure via an html page.

\section{Protocol of experimentation}
\label{sec:protocol_experimentation}

\subsection{Get one measure}
To be able to use these elements to get proper measure of the energy consumption of the GPU for different computations, we develloped programs and scripts to allow easy data gathering.

\begin{figure}[!ht]
\centering
\includegraphics[scale=0.4]{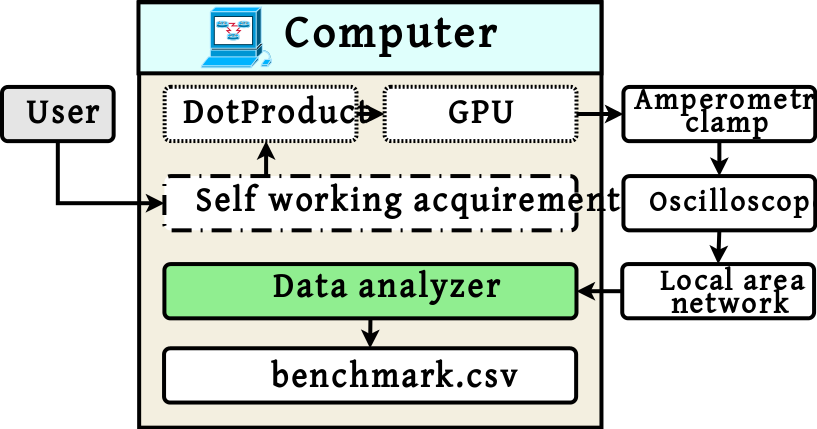}
\caption{Program process}
\label{fig:img:program_flow}
\end{figure}

The following schematic presents the flow to acquire one measure. The example program for a \verb+DotProduct+ kernel is given in Figure~\ref{fig:img:program_flow}.

The different steps to get these results are:
\begin{itemize}
\item The user launches the program that is using the DotProduct kernel.
\item The evolution of the current consumption of the GPU can be seen on the osciloscope screen.
\item When the execution finishes, the user can launch the dump of the data from the oscilloscope to the computer through the LAN.
\item The data acquired are a serie of point in Volt over Time. We have as many curves as amperometric clamps.
\item A data analyser script then combine these curves to compute the exact energy consumption of the GPU at each moment of the execution.
\end{itemize}

This method allows an easy acquisition of curves presenting the exact energy consumption at each moment of the program execution.
To get proper curves from one measure, we had to use simple signal processing. Indeed, our data where very noisy, mainly because we where working just above the lowest point of linearity of the amperometric clamps (the clamps we used where not made to measure small currents). These result where only used for visual rendering, this is why we only used an averaging filter to reduce the noise without deforming the curve too much.

\begin{figure}[!ht]
\centering
\includegraphics[scale=0.14]{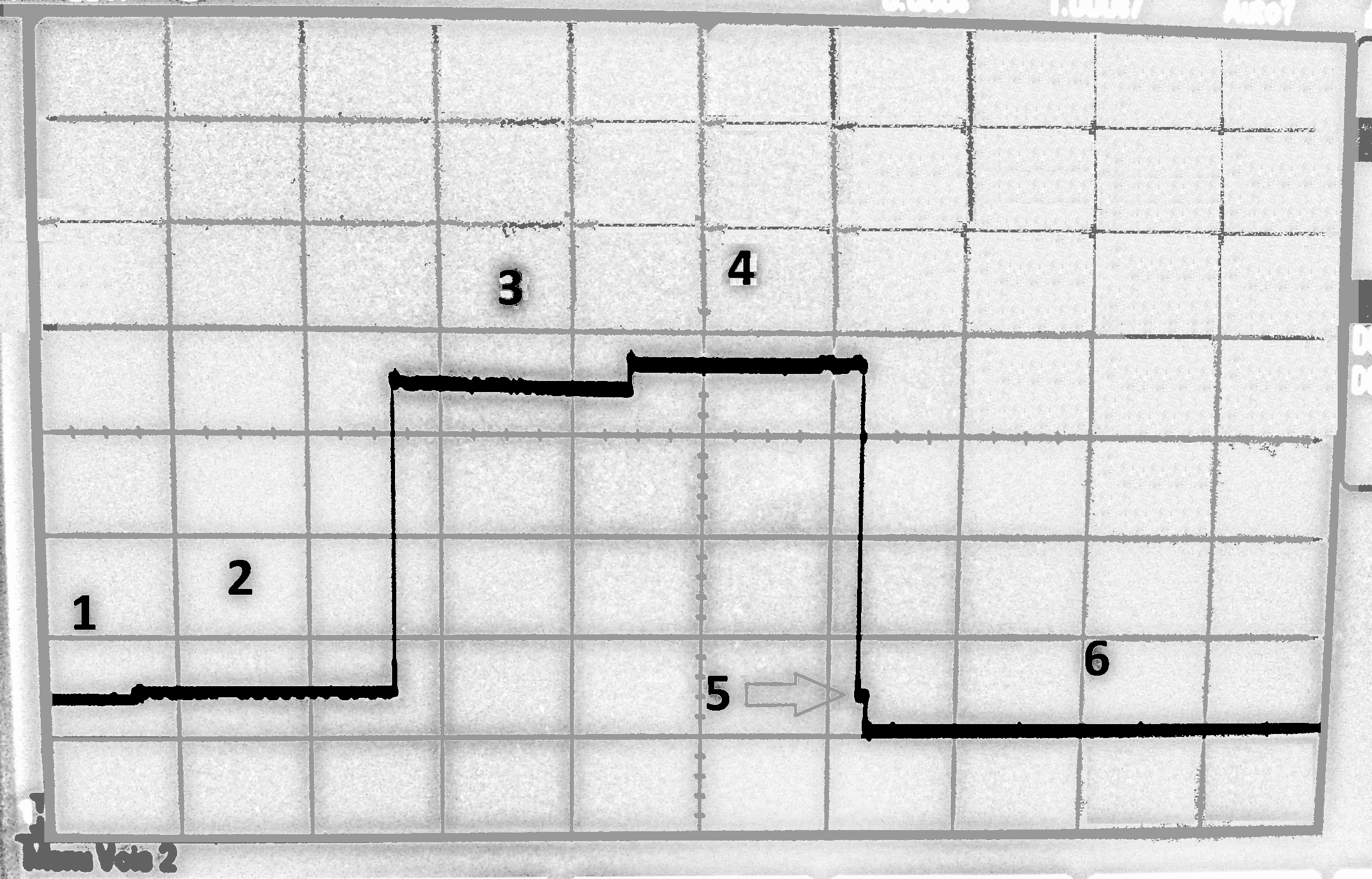}
\caption{Screen of the oscilloscope}
\label{fig:img:oscilloscope1}
\end{figure}

The kind of results that we got with this experiment are very interesting: above the view from the oscilloscope where we sum both signal from the 12V power-supply and use the low pass filter built in the osciloscope. A screen of the oscilloscope is drawn in~\ref{fig:img:oscilloscope1}.
\begin{enumerate}
\item The GPU is in idle mode. (CUDA not initialized yet)
\item We initialized some memory and copied elements (cudaMalloc and cudaMemCpy) then wait to be able to the the energy consumption.
\item Use of all the SM to perform additions.
\item Use of all the SM to perform multiplications.
\item Copy of the results back to the main memory.
\item The GPU is in idle mode. (program killed)
\end{enumerate}
The results obtained in Figure~\ref{fig:img:oscilloscope1} can be used to have a qualitative idea about the energy consumption but to get proper results, we need to get these values to the computer, perform several time the same computation to be sure the results remains the same and extract the uselful data in a more user-friendly way.

\subsection{Get a usefull set of data}

To be able to get lots of useful measure from this experimental setup, we needed to adress the three points presented above.

To get the data directly to the computer, we used an Ethernet-connected oscilloscope to be able to get and parse from the network the data we needed. The numerical oscilloscope we used also allowed us to directly compute the sum of several inputs before sending them over the network. In the setup we used, we limited the number of amperometric clamps to two: the power supplied by the 3.3V pins on the PCIe was not relevant compared to the 12V ones from the PCIe pins and the external allimentation. That way, we only used two entries from the oscilloscope that we sums before sending them to the computer. On the computer side we then multiply the value we got by 12 to get the energy used by the card.

\begin{figure}[!ht]
\centering
\includegraphics[scale=0.4]{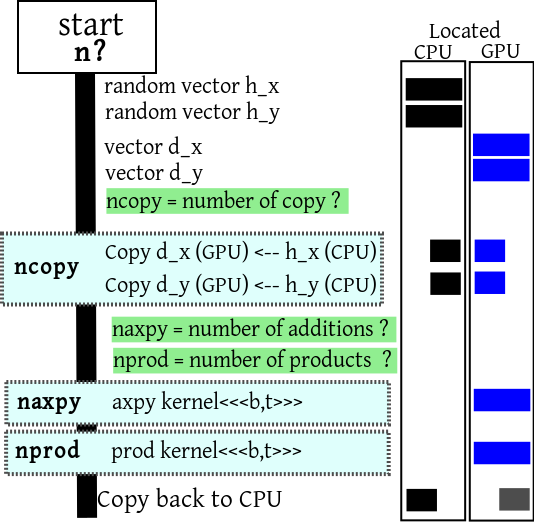}
\caption{Acquisition procedure}
\label{fig:img:algo_process_acquisition}
\end{figure}

We also wanted to be able to perform several times the same computation to be able to reduce the error of the results because the consumption usually changes by few watts between two measures. For that purpose, we created multiple scripts to automatize the acquisition:
\begin{itemize}
\item First the user need to specify the calculus that he wants, how many times each should run and the post processing function that should be use after each acquisition and after each serie of the same experiment.
\item The script launches each experiment following the user's specifications.
\item It launches the experimental CUDA program.
\item At the end of the execution, it launches the data gathering script to retrieve the data from the oscilloscope.
\item It then executes the selected post processing script.
\item After all the same experiments are finished, the corresponding post processing script is launched.
\item The script pause between two executions and ensure that the GPU charge is back to idle before starting a new loop.
\end{itemize}
The acquisition procedure is summarized in Figure~\ref{fig:img:algo_process_acquisition}.
Finally we wanted to make the result of our experiment more user friendly. To be able to do that, we used a property of the GPU energy consumption: it can be approximated with a step function. Figure~\ref{fig:img:measureexmpl} presents an example of data that we got from the oscilloscope:
\begin{figure}[!ht]
\centering
\includegraphics[scale=0.4]{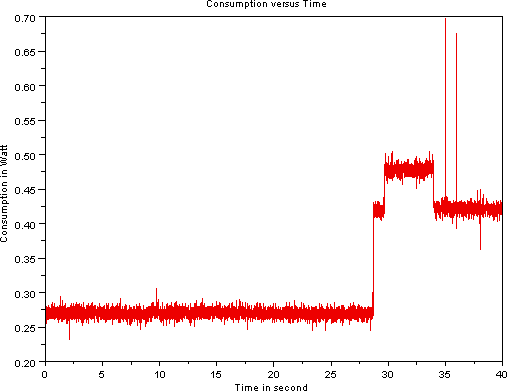}
\caption{Detection phase: energy consumption}
\label{fig:img:measureexmpl}
\end{figure}
An easy way to describe this graph is to have the time value for each step and the value between two jumps. We then used an algorithm that uses the variance variation to detect the step in the curves. To have good result even with noisy signals, we fed the algorithm with a sample of the noise from the sensor so it can take into account the standard deviation. This algorithm allows us to properly detect the steps from the acquired curves. An example of the output of this algorithm is given in Figure~\ref{fig:img:measureexmpl3}.
\begin{figure}[!ht]
\centering
\includegraphics[scale=0.40]{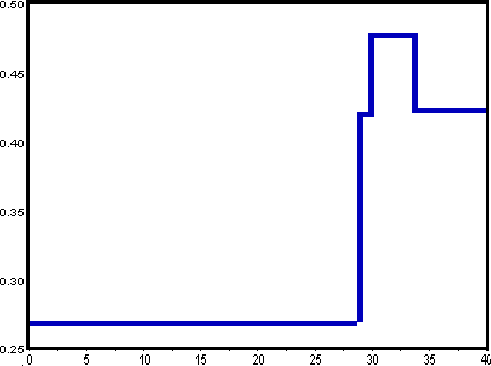}
\caption{Detection phase: approximate function of Figure~\ref{fig:img:measureexmpl}}
\label{fig:img:measureexmpl3}
\end{figure}

Each of our measurement is then described with a list of steps and of constant values. To be able to get a final value for one experiment, we then detected from all the measures of the same experiment, the ones with the most common pattern and then compute the mean of all the measures having this pattern (same number of steps and patterns for the values of each step).
We fully automatize this process with a configuration file describing all the experiments needed by the user. We could run as many experiments as needed, and summarize all the results in an excel spreadsheet containing for each experiment the power consumption for each phase.

\section{Numerical results}
\label{sec:exp_results}

To validate your protocol we propose to solve linear systems arising from finite element discretization of the gravity equations using the Conjugate Gradient (CG) method.
In this experiments we use CUDA library for GPU computations with real number arithmetics in double precision.

\subsection{Matrices tested}
\label{subsec:matrices_tested}

The sparse matrix we use arrises from a finite element method with $Q1$ discretisation of the Chicxulub crater in the Yucatan Peninsula. We modelised this 65 million year-old crater using the following Poisson equation
\begin{equation}
 \left\{
\begin{array}{l l l l l l}
-\Delta \Phi       & = & 4\pi G \delta\rho & = & f & \text{ on } \Omega\\
\mathcal{D}(\Phi) & = & b                  &   &   & \text{ on } \partial\Omega\\
\end{array} \right.\label{eqn:Poisson}
\end{equation}
where $\delta\rho$ is the density of the anomaly, $G$ the gravitational constant and $\mathcal{D}$ denotes Dirichlet boundary conditions. The computationnal domain $\Omega$ is a box of $200 \times 200 \times 10$ km$^3$ containing the crater of 10km deep and 180km large.\\

The properties of the set of matrices obtained from gravity model~\cite{cheikahamed:2013:inproceedings-3,ahamed2013stochastic} are described in Table~\ref{tab:sketches_matrices_gravimetry_tab} and an example of a common pattern is reported in Table~\ref{tab:sketches_matrices_gravimetry} where {\bf h} is the size of the matrix, {\bf nz} is the number of non-zero elements, {\bf density} is its density, i.e., the number of non-zero elements divided by the number of elements, {\bf bandwith} is the bandwith (upper and lower ones are equals for symmetric matrix), {\bf max\_row} is the maximum row density, {\bf nz/h} is the mean row density and {\bf nz/h stddev} is the standard deviation of the mean row density.
\begin{table}[!ht]  
\centering
\scriptsize
\renewcommand{\arraystretch}{1.3}
\renewcommand{\tabcolsep}{0.08cm}
\caption{Canvas of the matrices}
\label{tab:sketches_matrices_gravimetry_tab}
\begin{tabular}{ccccccc}
\hlinewd{1.0pt}
{\bf h} & {\bf nz} & {\bf density} & {\bf bandwidth} & {\bf max\_row} & {\bf nz/h} & {\bf nz/h stddev}\\
\hlinewd{1.0pt}
101168  & 2310912  & 0.022 & 6209  & 27 & 22.752 & 8.666 \\
296208  & 7101472  & 0.008 & 12869 & 27 & 23.975 & 7.544 \\
650848  & 16044032 & 0.004 & 21929 & 27 & 24.651 & 6.758 \\
848256  & 21073920 & 0.003 & 26225 & 27 & 24.844 & 6.504 \\
1213488 & 30434592 & 0.002 & 33389 & 27 & 25.080 & 6.171 \\
1325848 & 33321792 & 0.002 & 33389 & 27 & 25.132 & 6.094 \\
1587808 & 40077872 & 0.002 & 40019 & 27 & 25.241 & 5.930 \\
\hlinewd{1.0pt}
\end{tabular}
\end{table}

Table~\ref{tab:sketches_matrices_gravimetry} shows the sparse matrix pattern in the second column and in the third column are reported the properties of the matrix.
\begin{table}[!ht]  
\centering
\scriptsize
\caption{Example of gravity matrix}
\label{tab:sketches_matrices_gravimetry}
\begin{tabular}{m{5.8cm}m{6.2cm}}
\hlinewd{1.0pt}
\infosparse{img/png/ecpmashpc_gravi1325848}{gravi 1587808}{1587808}{40077872}{0.002}{40019}{27}{25.241}{5.930}\\
\hlinewd{1.0pt}
\end{tabular}
\end{table}

\subsection{Linear algebra operations}
\label{subsec:linear_algebra}

\begin{figure}[!ht]
\centering
\includegraphics[scale=0.23]{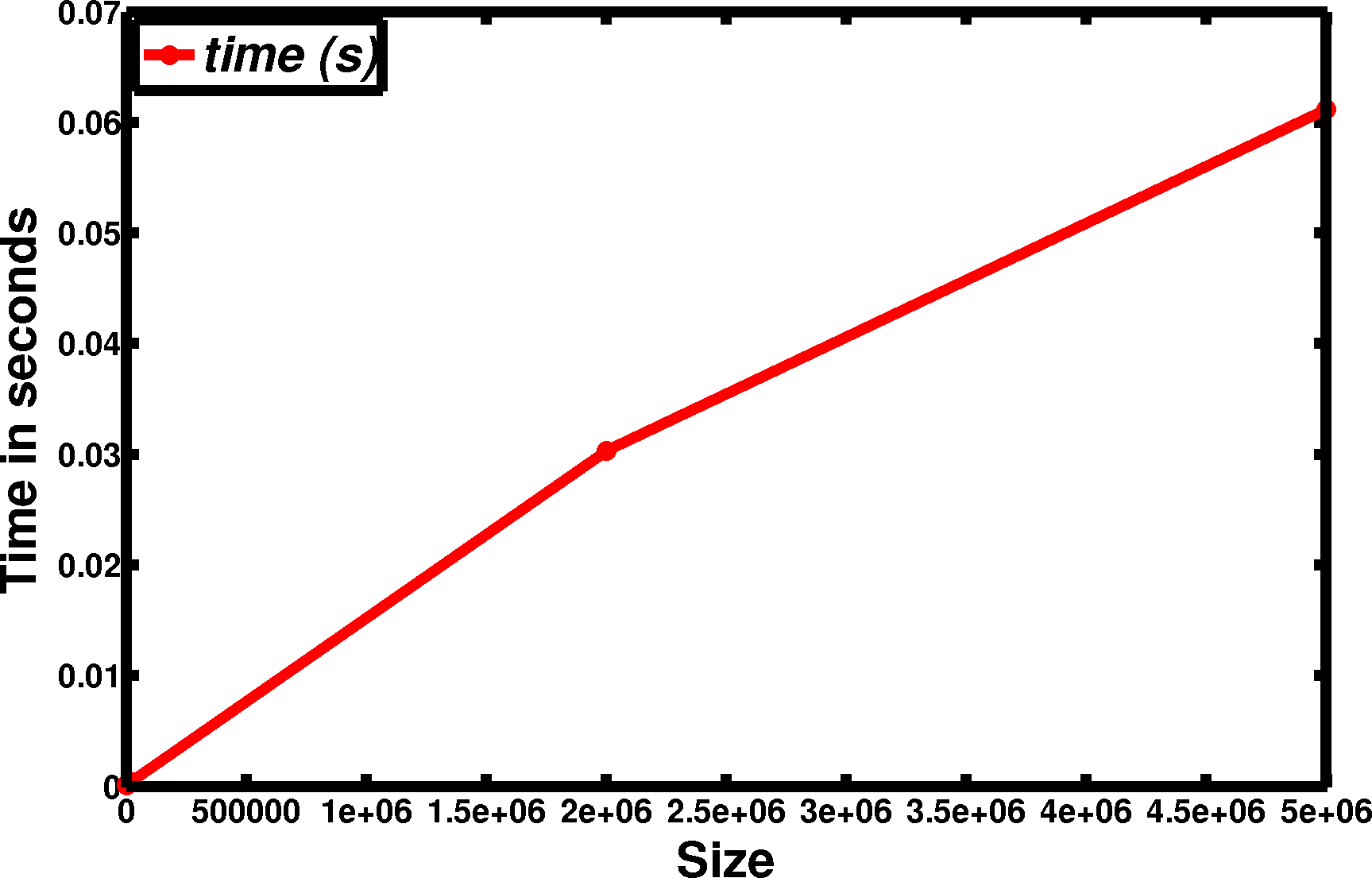}
\includegraphics[scale=0.23]{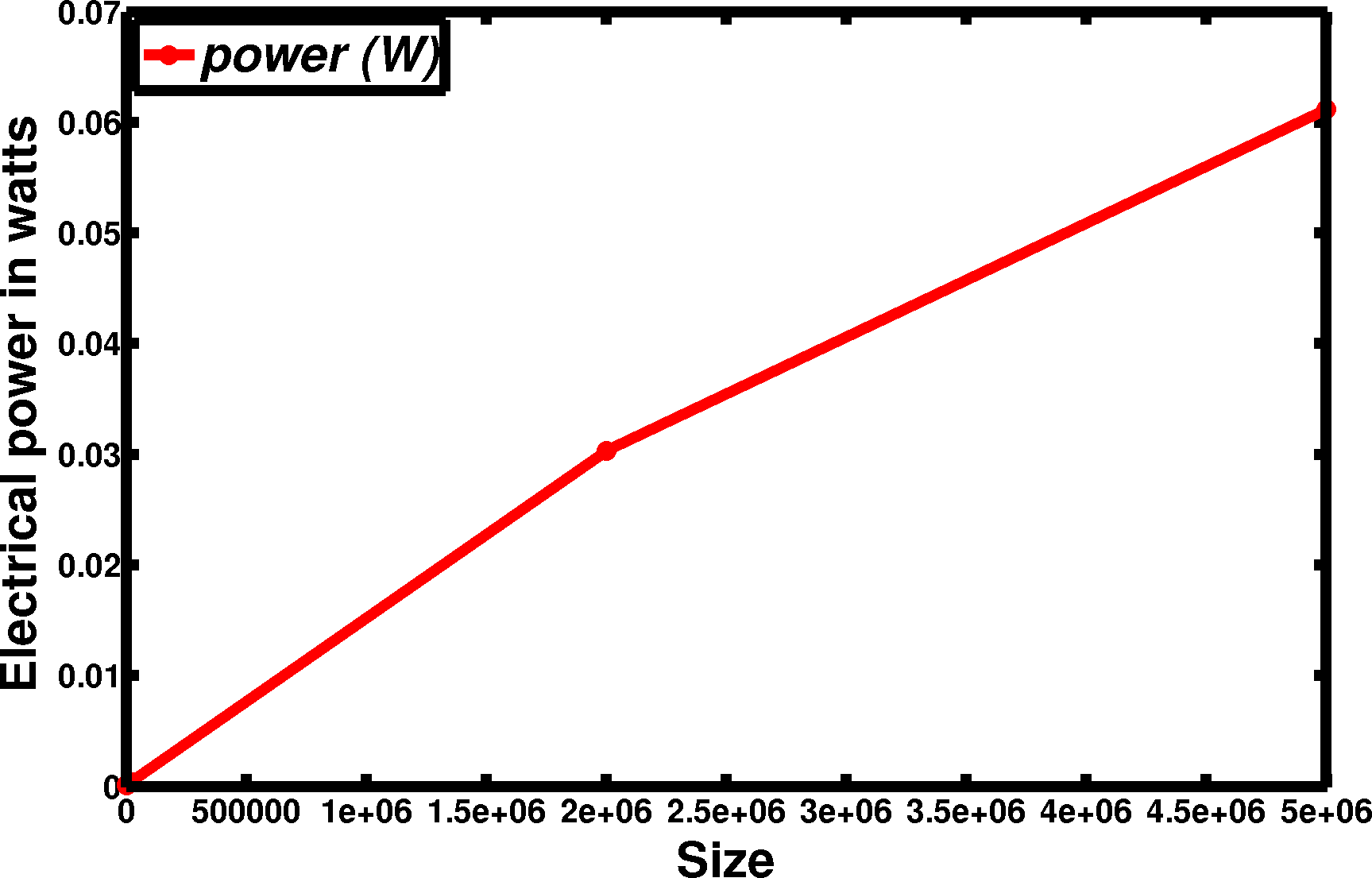}
\caption{Transfer data from CPU/GPU to GPU/CPU (Left: Time in seconds, Right: Power in Watt)}
\label{fig:img:time_power_data_transfer}
\end{figure}

Figure~\ref{fig:img:time_power_data_transfer} reports both the execution time in seconds (s) and power consumption in Watt (W) when data are transfered from CPU to GPU. 

Table~\ref{tab:daxpy_tab} and Table~\ref{tab:element_wise_product_tab} collect the running time in seconds (s) and power consumption (W) respectively for addition of vectors and element wise product operations.

\begin{table}[!ht]
\centering
\begin{tabular}{ccc}
\hline
{\bf size} & {\bf time (s)} & {\bf power (W)} \\
\hline
10,000 &    0.000006535 & 0.000420392 \\
2,000,000 & 0.000456454 & 0.027676176 \\
5,000,000 & 0.001142272 & 0.069889335 \\
\hline
\end{tabular}
\caption{Addition of vectors}
\label{tab:daxpy_tab}
\end{table}
\begin{table}[!ht]
\centering
\begin{tabular}{ccc}
\hline
{\bf size} & {\bf time (s)} & {\bf power (W)} \\
\hline
10,000   & 0.000006535 & 0.000541402 \\
2,000,000 & 0.000456454 & 0.046621519 \\
5,000,000 & 0.001142272 & 0.160956136 \\
\hline
\end{tabular}
\caption{Element wise product}
\label{tab:element_wise_product_tab}
\end{table}

\begin{figure}[!ht]
\centering
\includegraphics[scale=0.23]{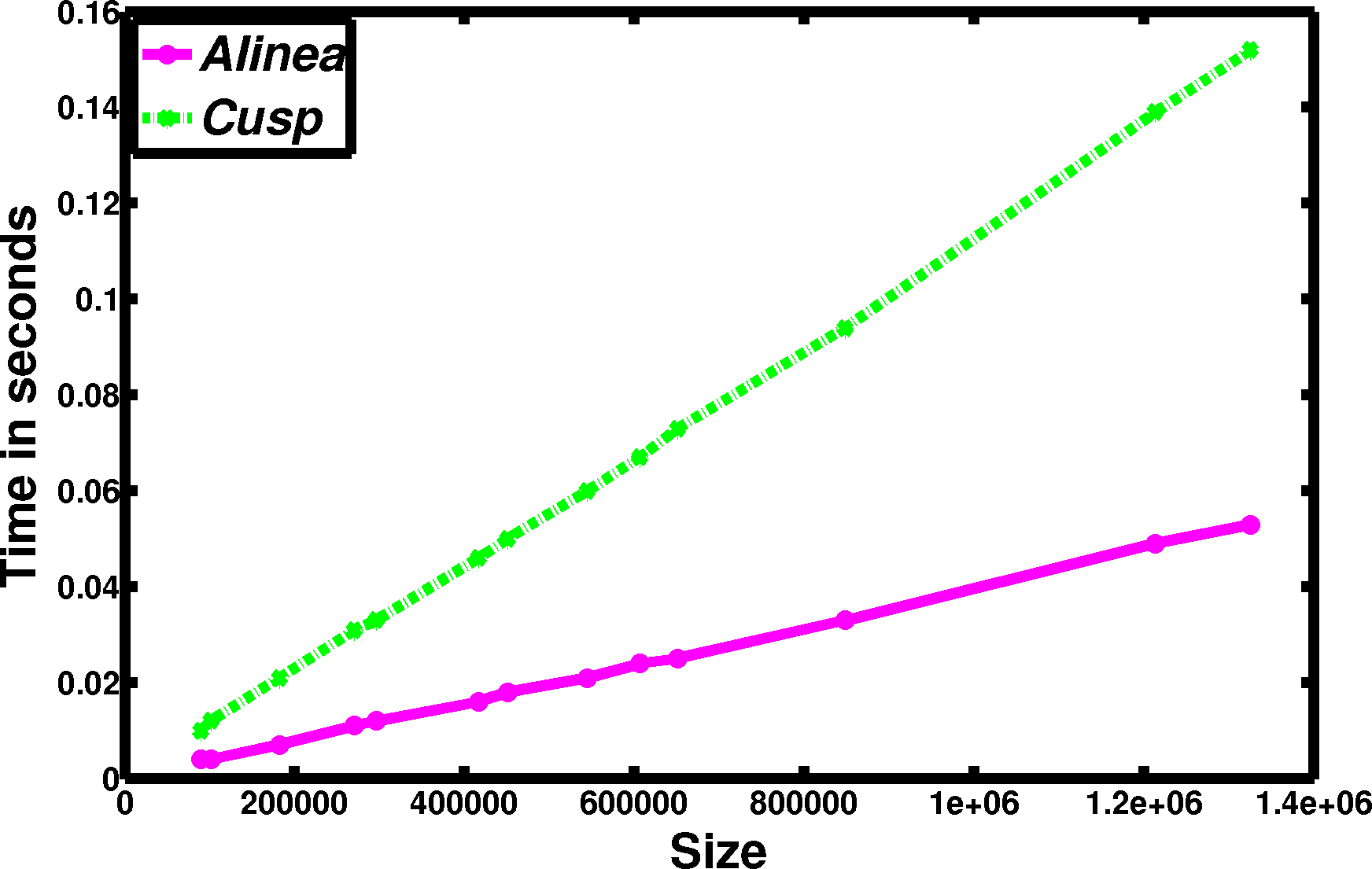}
\includegraphics[scale=0.23]{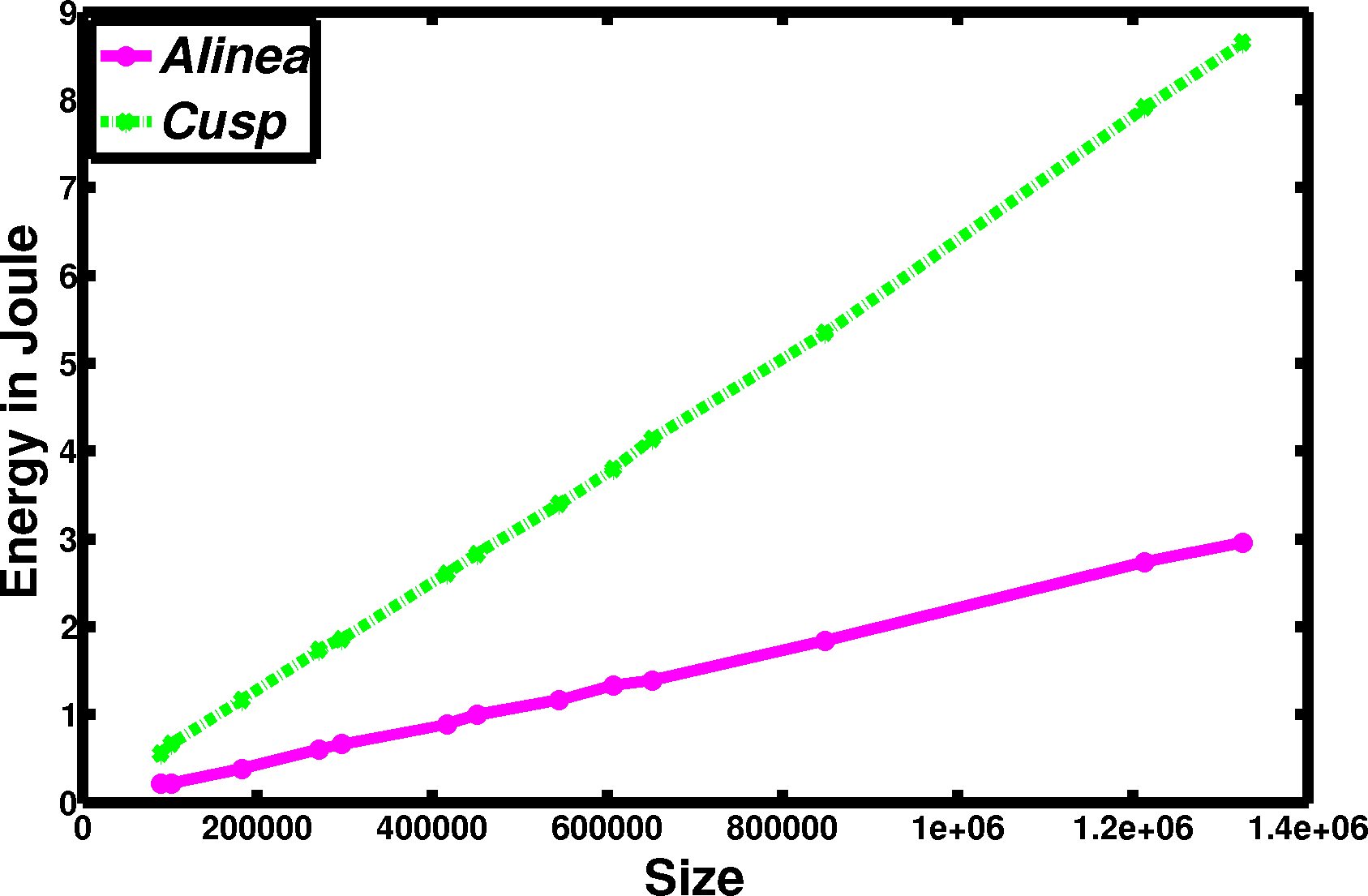}
\caption{Dot product (Left: Time in seconds, Right: Energy in Joule)}
\label{fig:img:time_energy_dot_product}
\end{figure}

In Table~\ref{tab:dot_product_tab} and Figure~\ref{fig:img:time_energy_dot_product} are reported the execution time in seconds, power consumption in Watt (W) and energy in Joule (J) for the dot product operations. We compare the results obtained by Alinea, our research library described in~\cite{cheikahamed:2012:inproceedings-1,cheikahamed:2012:inproceedings-2,cheikahamed:2013:inproceedings-3,cheikahamed:2013:inproceedings-4} against the Cusp library~\cite{GPU:CUSP0.3.0:2012}. As we can see in this Table, Alinea outperforms Cusp in terms of energy consumption.
\begin{table}[!ht]
\centering
\begin{tabular}{ccccccc}
\hline
\multicolumn{1}{c|}{}    & \multicolumn{3}{c|}{ {\bf Alinea} } & \multicolumn{3}{c}{ {\bf Cusp } } \\
\multicolumn{1}{c|}{\bf n} & {\bf time (s)} & {\bf P (W)} & \multicolumn{1}{c|}{\bf E (J)} & {\bf time (s)} & {\bf P (W)} & {\bf E (J)} \\
\hline
101,168 & 0.004 & 54.672 & 0.227 & 0.012 & 55.717 & 0.659 \\
296,208 & 0.012 & 55.427 & 0.645 & 0.033 & 56.457 & 1.875 \\
650,848 & 0.025 & 55.604 & 1.418 & 0.073 & 56.763 & 4.130 \\
848,256 & 0.033 & 55.865 & 1.842 & 0.094 & 56.884 & 5.351 \\
1,213,488 & 0.049 & 55.884 & 2.718 & 0.139 & 56.927 & 7.901 \\
1,325,848 & 0.053 & 55.863 & 2.968 & 0.152 & 56.906 & 8.628 \\
1587808 & 0.064 & 55.892 & 3.555 & 0.182 & 56.935 & 10.335 \\
\hline
\end{tabular}
\caption{Dot product}
\label{tab:dot_product_tab}
\end{table}

\begin{figure}[!ht]
\centering
\includegraphics[scale=0.23]{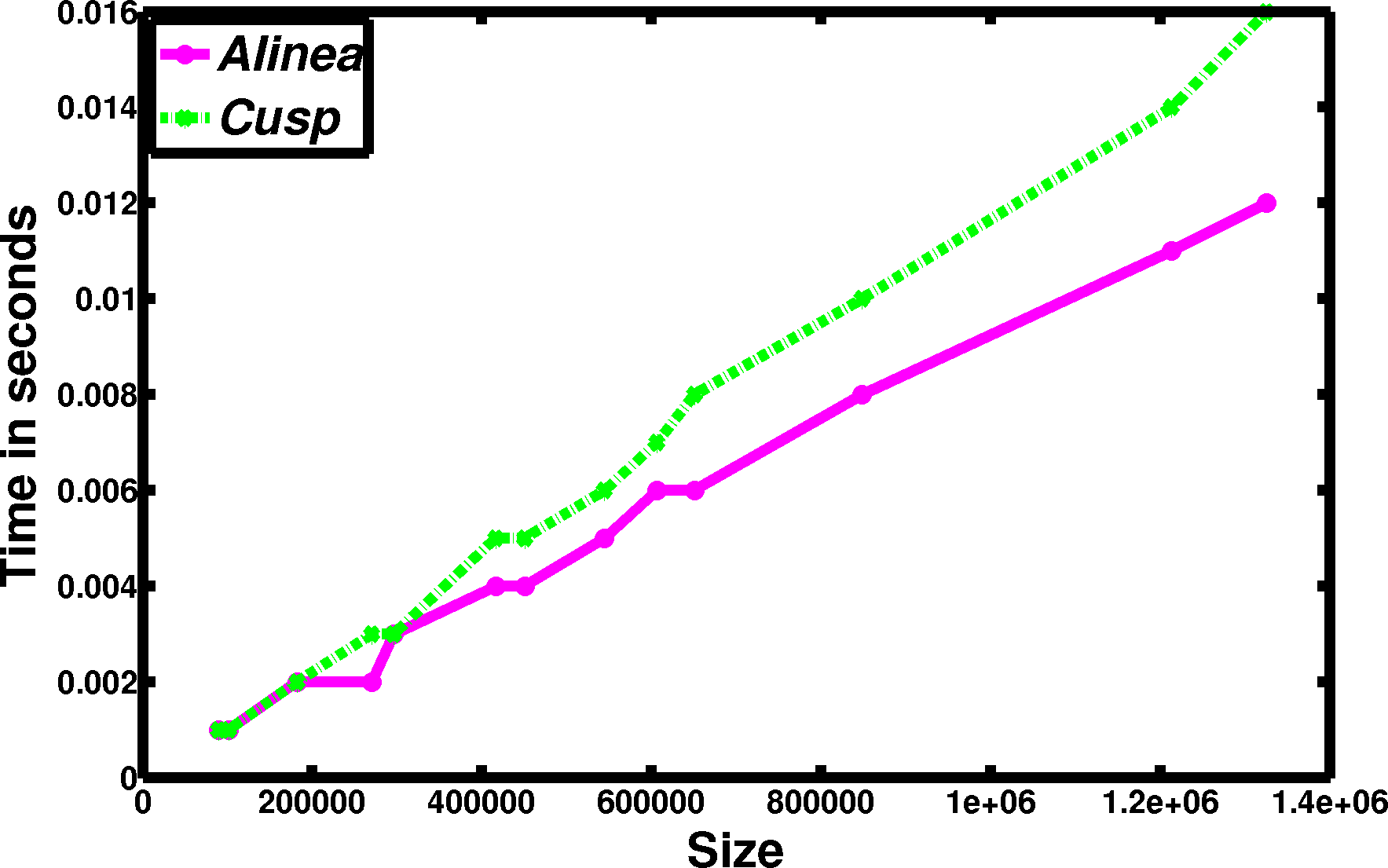}
\includegraphics[scale=0.23]{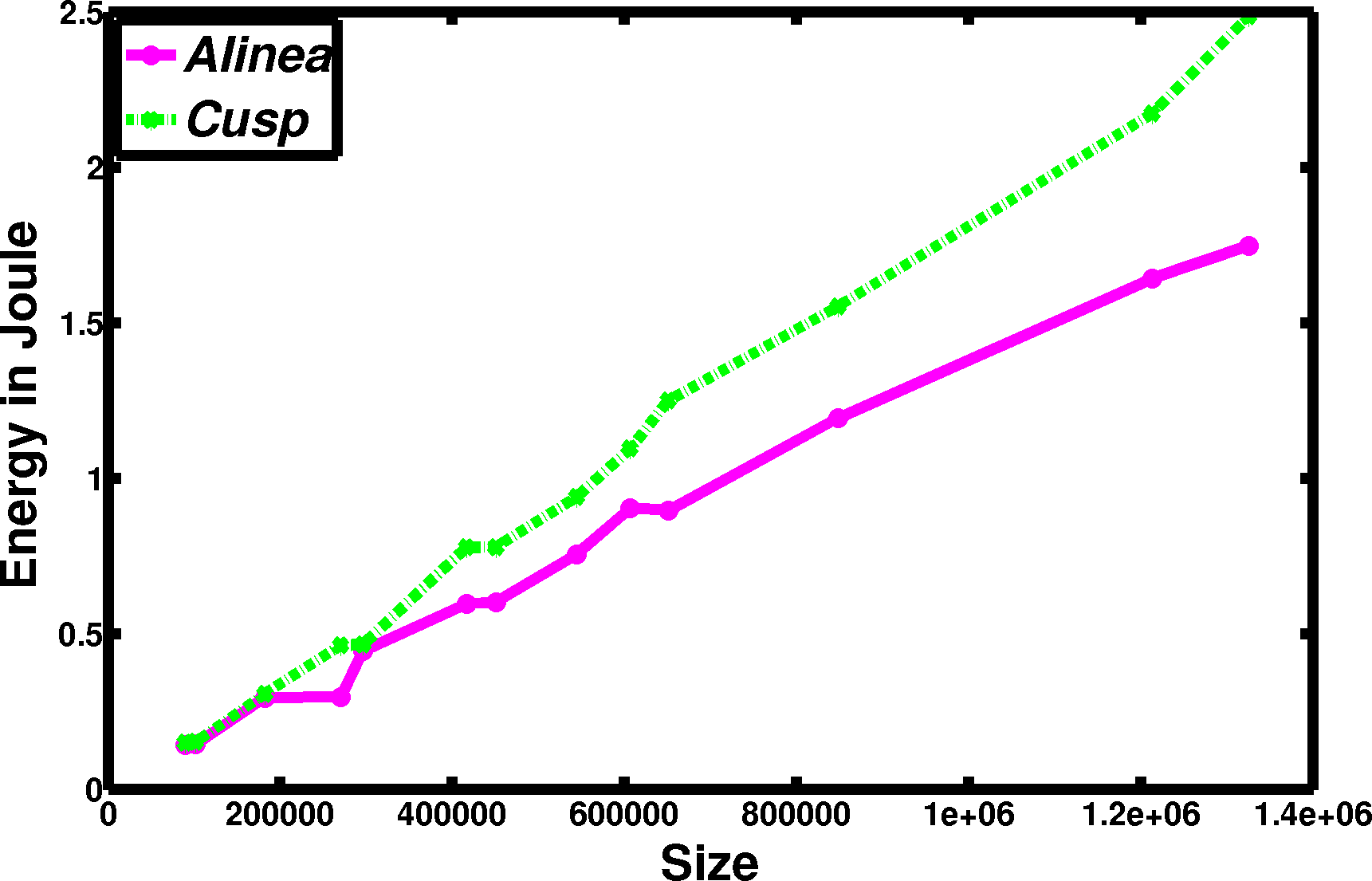}
\caption{Sparse matrix-vector product (Left: Time in seconds, Right: Energy in Joule)}
\label{fig:img:time_energy_spmv}
\end{figure}

Table~\ref{tab:spmv_tab} and Figure~\ref{fig:img:time_energy_spmv} present the sparse matrix-vector product (SpMV) execution in seconds, power (W) and energy (J) consumption.  The SpMV is the most consuming operations in terms of computing time~\cite{PARALLEL:BFFGL:2009,GPU:CSVGM:2014,GPU:BG:2008,GPU:BG:2009}.
\begin{table}[!ht]
\centering
\begin{tabular}{ccccccc}
\hline
\multicolumn{1}{c|}{}    & \multicolumn{3}{c|}{ {\bf Alinea} } & \multicolumn{3}{c}{ {\bf Cusp } } \\
\multicolumn{1}{c|}{\bf n} & {\bf time (s)} & {\bf P (W)} & \multicolumn{1}{c|}{\bf E (J)} & {\bf time (s)} & {\bf P (W)} & {\bf E (J)} \\
\hline
101,168 & 0.001 & 144.229 & 0.136 & 0.001 & 151.423 & 0.183 \\
296,208 & 0.003 & 148.591 & 0.401 & 0.003 & 154.756 & 0.541 \\
650,848 & 0.006 & 149.420 & 0.891 & 0.008 & 156.229 & 1.198 \\
848,256 & 0.008 & 149.289 & 1.176 & 0.010 & 155.357 & 1.552 \\
1,213,488 & 0.011 & 149.354 & 1.676 & 0.014 & 155.280 & 2.219 \\
1,325,848 & 0.012 & 145.663 & 1.807 & 0.016 & 155.419 & 2.427 \\
1587808 & 0.015 & 147.669 & 2.165 & 0.019 & 155.728 & 2.913 \\
\hline
\end{tabular}
\caption{SpMV with CSR format}
\label{tab:spmv_tab}
\end{table}
Table~\ref{tab:spmv_tab} and Figure~\ref{fig:img:time_energy_spmv} confirm the robustness and efficiency of our library in terms of computing time and energy consumption for double precision arithmetics.

\subsection{Conjugate gradient method}
\label{subsec:conj_gradient}

Kwowing the efficiency of linear algebra operations, we now analysis the obtained results with the conjugate gradient, which an efficient iterative Krylov~\cite{saad_iterative_2003,GPU:BCK:2011,PARALLEL:AHR:2010} method for solving large sparse  linear systems~\cite{GPU:LS:2010,PARALLEL:AHR:2010,gaikwad_parallel_2010} with positive definite matrices.
The obtained results for classical CG are collected in Table~\ref{tab:iterative_krylov_cg_tab} and Figure~\ref{tab:iterative_krylov_cg_tab}.
\begin{table}[!ht]
\centering
\begin{tabular}{cccccccc}
\hline
\multicolumn{1}{c|}{} & \multicolumn{1}{c|}{} & \multicolumn{3}{c|}{ {\bf Alinea} } & \multicolumn{3}{c}{ {\bf Cusp } } \\
\multicolumn{1}{c|}{n} & \multicolumn{1}{c|}{\bf \#iter} & {\bf time (s)} & {\bf P (W)} & \multicolumn{1}{c|}{\bf E (J)} & {\bf time (s)} & {\bf P (W)} & {\bf E (J)} \\
\hline
101,168 & 203 & 0.313 & 100.201 & 31.363 & 0.471 & 110.271 & 51.888 \\
296,208 & 299 & 1.155 & 123.190 & 142.284 & 1.436 & 133.129 & 191.199 \\
650,848 & 394 & 2.460 & 133.772 & 329.078 & 3.741 & 141.852 & 530.609 \\
848,256 & 435 & 4.336 & 137.843 & 597.618 & 5.250 & 145.828 & 765.595 \\
1,213,488 & 491 & 6.872 & 141.169 & 970.115 & 8.396 & 146.889 & 1,233.248 \\
1,325,848 & 538 & 8.180 & 141.747 & 1,159.492 & 9.945 & 150.403 & 1,495.759 \\
\hline
\end{tabular}
\caption{CG with CSR format}
\label{tab:iterative_krylov_cg_tab}
\end{table}
\begin{figure}[!ht]
\centering
\includegraphics[scale=0.23]{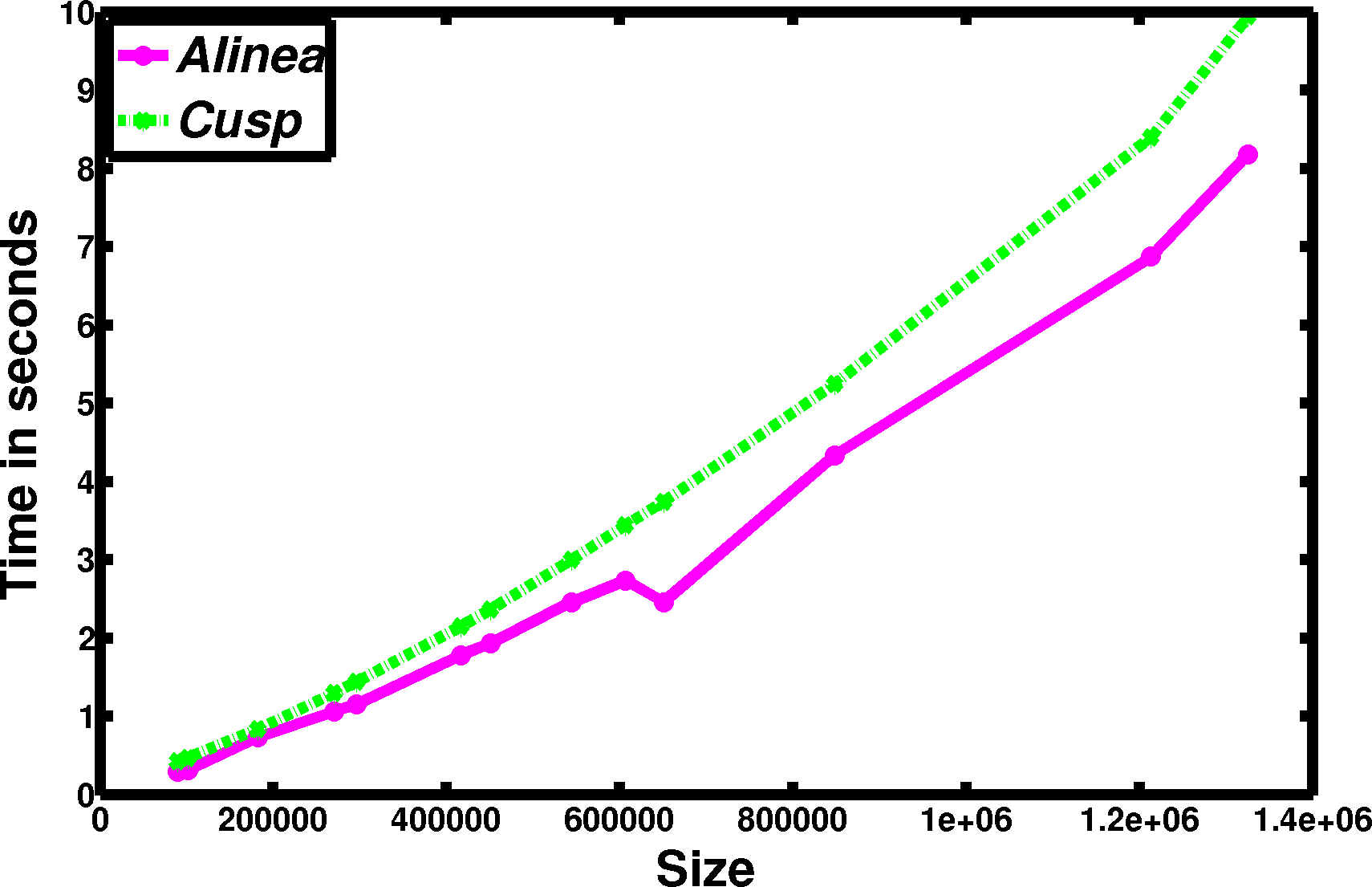}
\includegraphics[scale=0.23]{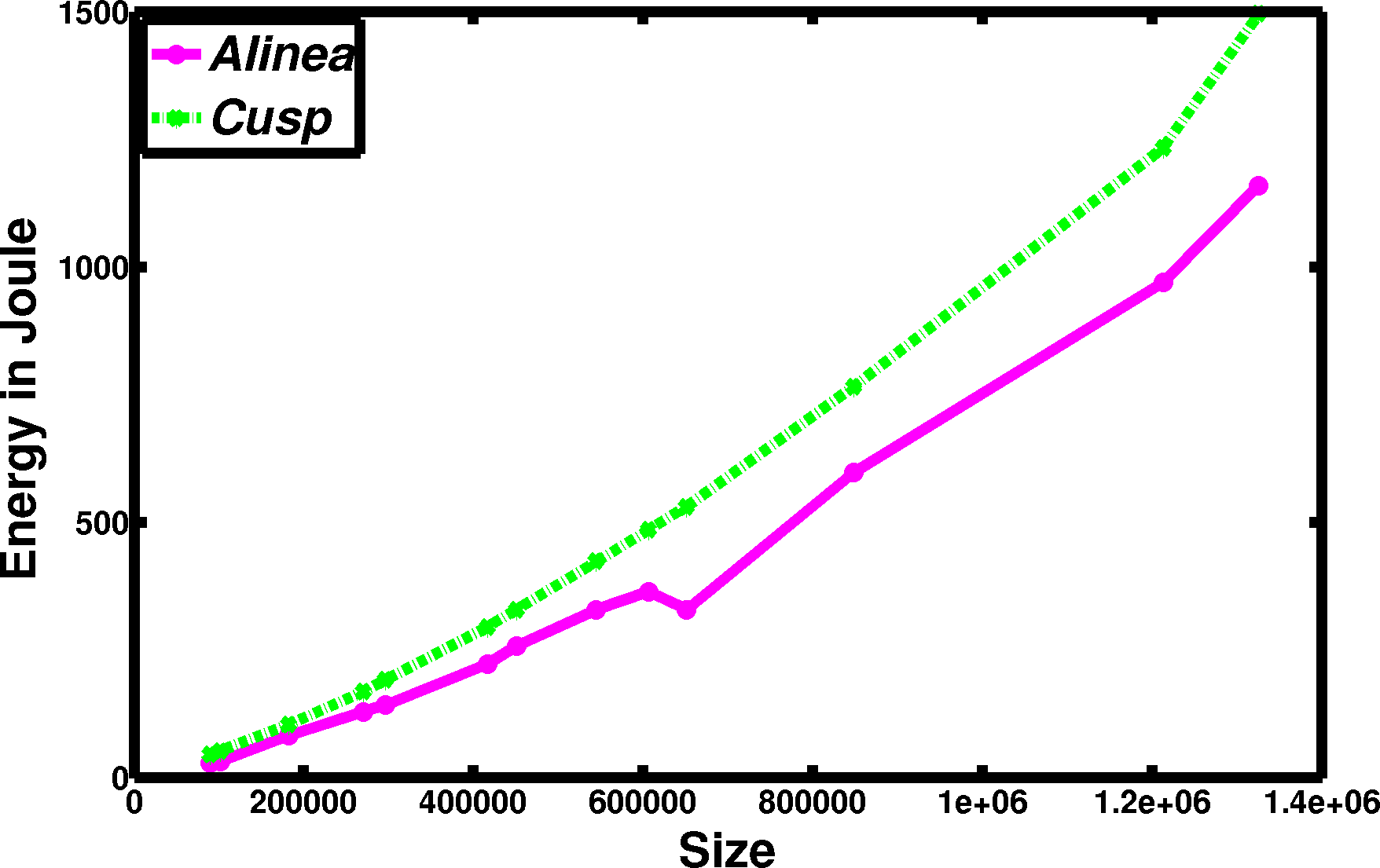}
\caption{Conjugate gradient (Left: Time in seconds, Right: Energy in Joule)}
\label{fig:img:time_energy_solver}
\end{figure}
As we can see in Table~\ref{tab:iterative_krylov_cg_tab} and Figure~\ref{tab:iterative_krylov_cg_tab}, Alinea is definitely better than Cusp in terms of energy consumption: Alinea is about 30\% more efficient than Cusp.

In the previous analysis no preconditioner have been used.
Anyway, it is fundamental to mention that CG can be used very efficiently on GPU with some preconditioning based on domain decomposition methods~\cite{magoules:journal-auth:21,magoules:journal-auth:16}.
For this purpose, special tuned interface conditions between the subdomains are usually defined like in the FETI method~\cite{magoules:journal-auth:24} or like in the Schwarz method with homogeneous coefficients~\cite{magoules:journal-auth:10,magoules:journal-auth:9} or with heterogeneous coefficients~\cite{magoules:journal-auth:23,magoules:journal-auth:18,magoules:journal-auth:14}.
Othe alternative to tune the interface conditions exists based on patch substructuring as introduced for acoustics in~\cite{magoules:proceedings-auth:6,magoules:journal-auth:29,magoules:journal-auth:12}
or for linear elasticity in~\cite{magoules:journal-auth:20,magoules:journal-auth:17}.
The associated problem, with the Lagrange multipliers, condensed on the interface between the subdomains, is then solved with an interative algorithm on the CPU.
At each iteration, each subproblem defined in each subdomain is solved with the CG method on the GPU.
This CPU-GPU implementation of domain decomposition methods were first proposed for the FETI method in~\cite{Papadrakakis:2011} and for the optimized Schwarz method in~\cite{cheikahamed:2013:inproceedings-4}.

\section{Concluding remarks}
\label{sec:concluding_remarks}

In this paper, a novel analysis allows us to compare algorithms with a new key variable: their exact energy consumption. With the new evolution of high performance computing and the increasing need to reduce the energy consumption, this new possibility to know the expected energy consumption of an algorithm, when developping it, is very useful.
Moreover this analysis can easily be generalized and used to predict the energy consumption of any algorithms when the operations involved within this algorithm are known.
This is an incredibly powerful tool when designing an algorithm to be able to predict its execution time and its energy consumption!

\section*{Acknowledgment}
The authors acknowledge the CUDA Research Center at Ecole Centrale Paris (France) for its support and for providing the computing facilities.

\bibliography{bib/hpcc2014_energy_ad_ac_fm,bib/MAGOULES-JOURNAL1,bib/MAGOULES-PROCEEDINGS1}

\begin{thebibliography}{10}

\bibitem{PARALLEL:AHR:2010}
H.~Anzt, V.~Heuveline, and B.~Rocker.
\newblock Mixed precision iterative refinement methods for linear systems:
  Convergence analysis based on {K}rylov subspace methods.
\newblock In K.~J\'onasson, editor, {\em PARA (2)}, volume 7134 of {\em Lecture
  Notes in Computer Science}, pages 237--247. Springer, 2010.

\bibitem{GPU:BCK:2011}
J.~M. Bahi, R.~Couturier, and L.~Z. Khodja.
\newblock Parallel {GMRES} implementation for solving sparse linear systems on
  {GPU} clusters.
\newblock In {\em Proceedings of the 19th High Performance Computing Symposia},
  pages 12--19, San Diego, CA, USA, 2011. Society for Computer Simulation
  International.

\bibitem{GPU:BYFWA:2009}
A.~Bakhoda, G.~Yuan, W.~Fung, H.~Wong, and T.~Aamodt.
\newblock Analyzing {CUDA} workloads using a detailed {GPU} simulator.
\newblock In {\em Performance Analysis of Systems and Software, 2009. ISPASS
  2009. IEEE International Symposium on}, pages 163--174, April 2009.

\bibitem{GPU:BG:2008}
N.~Bell and M.~Garland.
\newblock Efficient sparse matrix-vector multiplication on {CUDA}.
\newblock Nvidia Technical Report NVR-2008-004, Nvidia Corporation, 2008.

\bibitem{GPU:BG:2009}
N.~Bell and M.~Garland.
\newblock Implementing sparse matrix-vector multiplication on
  throughput-oriented processors.
\newblock In {\em Proceedings of the Conference on High Performance Computing
  Networking, Storage and Analysis (SC'09)}, pages 1--11, New York, NY, USA,
  2009. ACM.

\bibitem{GPU:CUSP0.3.0:2012}
N.~Bell and M.~Garland.
\newblock Cusp: Generic parallel algorithms for sparse matrix and graph
  computations, 2012.
\newblock Available on line at: \url{http://cusplibrary.github.io/} (accessed
  on \today).

\bibitem{Budruk:2003:PES:861280}
R.~Budruk, D.~Anderson, and E.~Solari.
\newblock {\em PCI Express System Architecture}.
\newblock Pearson Education, 2003.

\bibitem{PARALLEL:BFFGL:2009}
A.~Buluç, J.~T. Fineman, M.~Frigo, J.~R. Gilbert, and C.~E. Leiserson.
\newblock Parallel sparse matrix-vector and matrix-transpose-vector
  multiplication using compressed sparse blocks.
\newblock In F.~M. auf~der Heide and M.~A. Bender, editors, {\em SPAA}, pages
  233--244. ACM, 2009.

\bibitem{GPU:CSVGM:2014}
A.~F. Camargos, V.~C. Silva, J.-M. Guichon, and G.~Meunier.
\newblock Iterative solution on {GPU} of linear systems arising from the a-v
  edge-fea of time-harmonic electromagnetic phenomena.
\newblock In {\em Parallel, Distributed and Network-Based Processing (PDP),
  2014 22nd Euromicro International Conference on}, pages 365--371, Feb 2014.

\bibitem{cheikahamed:2012:inproceedings-2}
A.-K. Cheik~Ahamed and F.~Magoul\`es.
\newblock Fast sparse matrix-vector multiplication on graphics processing unit
  for finite element analysis.
\newblock In {\em High Performance Computing and Communication 2012 IEEE 9th
  International Conference on Embedded Software and Systems (HPCC-ICESS), 2012
  IEEE 14th International Conference on}, pages 1307--1314. IEEE Computer
  Society, 2012.

\bibitem{cheikahamed:2012:inproceedings-1}
A.-K. Cheik~Ahamed and F.~Magoul\`es.
\newblock Iterative methods for sparse linear systems on graphics processing
  unit.
\newblock In {\em High Performance Computing and Communication 2012 IEEE 9th
  International Conference on Embedded Software and Systems (HPCC-ICESS), 2012
  IEEE 14th International Conference on}, pages 836--842. IEEE Computer
  Society, june 2012.

\bibitem{cheikahamed:2013:inproceedings-3}
A.-K. Cheik~Ahamed and F.~Magoul\`es.
\newblock Iterative {K}rylov methods for gravity problems on graphics
  processing unit.
\newblock In {\em Distributed Computing and Applications to Business,
  Engineering Science (DCABES), 2013 12th International Symposium on}, pages
  16--20. IEEE Computer Society, 2013.

\bibitem{cheikahamed:2013:inproceedings-4}
A.-K. Cheik~Ahamed and F.~Magoul\`es.
\newblock Schwarz method with two-sided transmission conditions for the gravity
  equations on graphics processing unit.
\newblock In {\em Distributed Computing and Applications to Business,
  Engineering Science (DCABES), 2013 12th International Symposium on}, pages
  105--109. IEEE Computer Society, 2013.

\bibitem{ahamed2013stochastic}
A.-K. Cheik~Ahamed and F.~Magoul\`es.
\newblock A stochastic-based optimized {S}chwarz method for the gravimetry
  equations on {GPU} clusters.
\newblock In {\em Domain Decomposition Methods in Science and Engineering XXI}.
  Springer, 2014.

\bibitem{ref:cheikahamed:2012:inproceedings-1:DARG:2013}
L.~Djinevski, S.~Arsenovski, S.~Ristov, and M.~Gusev.
\newblock Performance drawbacks for matrix multiplication using set associative
  cache in {GPU} devices.
\newblock In {\em Information \& Communication Technology Electronics \&
  Microelectronics (MIPRO), 2013 36th International Convention on}, pages
  193--198. IEEE, 2013.

\bibitem{gaikwad_parallel_2010}
A.~Gaikwad and I.~Toke.
\newblock Parallel iterative linear solvers on {GPU:} a financial engineering
  case.
\newblock In {\em 2010 18th Euromicro International Conference on Parallel,
  Distributed and Network-Based Processing ({PDP)}}, pages 607--614, Feb. 2010.

\bibitem{GPU:green500:2013}
Green500, 2013.
\newblock Available on line at: \url{http://www.green500.org} (accessed on
  \today).

\bibitem{GPU:LS:2010}
R.~Li and Y.~Saad.
\newblock {GPU}-accelerated preconditioned iterative linear solvers.
\newblock Technical Report umsi-2010-112, Minnesota Supercomputer Institute,
  University of Minnesota, Minneapolis, MN, 2010.

\bibitem{LDHDBX2011}
W.~Liu, Z.~Du, Y.~Hiao, B.~David~A., and C.~Xu.
\newblock A waterfall model to achieve energy efficient tasks mapping for large
  scale {GPU} clusters.
\newblock In {\em Parallel and Distributed Processing Workshops and Phd Forum
  (IPDPSW), 2011 IEEE International Symposium on}, pages 82 -- 92, May 2011.

\bibitem{magoules:journal-auth:14}
Y.~Maday and F.~Magoul\`es.
\newblock Non-overlapping additive {S}chwarz methods tuned to highly
  heterogeneous media.
\newblock {\em Comptes Rendus {\`a} l'Acad{\'e}mie des Sciences},
  341(11):701--705, 2005.

\bibitem{magoules:journal-auth:16}
Y.~Maday and F.~Magoul\`es.
\newblock Absorbing interface conditions for domain decomposition methods: a
  general presentation.
\newblock {\em Computer Methods in Applied Mechanics and Engineering},
  195(29--32):3880--3900, 2006.

\bibitem{magoules:journal-auth:18}
Y.~Maday and F.~Magoul\`es.
\newblock Improved ad hoc interface conditions for {S}chwarz solution procedure
  tuned to highly heterogeneous media.
\newblock {\em Applied Mathematical Modelling}, 30(8):731--743, 2006.

\bibitem{magoules:journal-auth:24}
Y.~Maday and F.~Magoul\`es.
\newblock Optimal convergence properties of the {FETI} domain decomposition
  method.
\newblock {\em International Journal for Numerical Methods in Fluids},
  55(1):1--14, 2007.

\bibitem{magoules:journal-auth:23}
Y.~Maday and F.~Magoul\`es.
\newblock Optimized {S}chwarz methods without overlap for highly heterogeneous
  media.
\newblock {\em Computer Methods in Applied Mechanics and Engineering},
  196(8):1541--1553, 2007.

\bibitem{magoules:journal-auth:9}
F.~Magoul\`es, P.~Iv\'anyi, and B.~Topping.
\newblock Convergence analysis of {S}chwarz methods without overlap for the
  {H}elmholtz equation.
\newblock {\em Computers and Structures}, 82(22):1835--1847, 2004.

\bibitem{magoules:journal-auth:10}
F.~Magoul\`es, P.~Iv\'anyi, and B.~Topping.
\newblock Non-overlapping {S}chwarz methods with optimized transmission
  conditions for the {H}elmholtz equation.
\newblock {\em Computer Methods in Applied Mechanics and Engineering},
  193(45--47):4797--4818, 2004.

\bibitem{magoules:journal-auth:21}
F.~Magoul\`es and F.-X. Roux.
\newblock Lagrangian formulation of domain decomposition methods: a unified
  theory.
\newblock {\em Applied Mathematical Modelling}, 30(7):593--615, 2006.

\bibitem{magoules:journal-auth:12}
F.~Magoul\`es, F.-X. Roux, and L.~Series.
\newblock Algebraic way to derive absorbing boundary conditions for the
  {H}elmholtz equation.
\newblock {\em Journal of Computational Acoustics}, 13(3):433--454, 2005.

\bibitem{magoules:journal-auth:17}
F.~Magoul\`es, F.-X. Roux, and L.~Series.
\newblock Algebraic approximation of {D}irichlet-to-{N}eumann maps for the
  equations of linear elasticity.
\newblock {\em Computer Methods in Applied Mechanics and Engineering},
  195(29--32):3742--3759, 2006.

\bibitem{magoules:journal-auth:20}
F.~Magoul\`es, F.-X. Roux, and L.~Series.
\newblock Algebraic {D}irichlet-to-{N}eumann mapping for linear elasticity
  problems with extreme contrasts in the coefficients.
\newblock {\em Applied Mathematical Modelling}, 30(8):702--713, 2006.

\bibitem{magoules:journal-auth:29}
F.~Magoul\`es, F.-X. Roux, and L.~Series.
\newblock Algebraic approach to absorbing boundary conditions for the
  {H}elmholtz equation.
\newblock {\em International Journal of Computer Mathematics}, 84(2):231--240,
  2007.

\bibitem{GPU:CUDA4.0:2011}
{Nvidia Corporation}.
\newblock {\em CUDA Toolkit Reference MANUAL}, 4.0 edition, 2011.
\newblock Available on line at:
  \url{http://developer.nvidia.com/cuda-toolkit-40} (accessed on \today).

\bibitem{Papadrakakis:2011}
M.~Papadrakakis, G.~Stavroulakis, and A.~Karatarakis.
\newblock {A new era in scientific computing: Domain decomposition methods in
  hybrid CPU–GPU architectures}.
\newblock {\em Computer Methods in Applied Mechanics and Engineering},
  200(13-16):1490--1508, Mar. 2011.

\bibitem{magoules:proceedings-auth:6}
F.-X. Roux, F.~Magoul\`es, L.~Series, and Y.~Boubendir.
\newblock Approximation of optimal interface boundary conditions for
  two-{L}agrange multiplier {FETI} method.
\newblock In R.~Kornhuber, R.~Hoppe, J.~P\'eriaux, O.~Pironneau, O.~Widlund,
  and J.~Xu, editors, {\em Proceedings of the 15th International Conference on
  Domain Decomposition Methods, Berlin, Germany, July 21-15, 2003}, {Lecture
  Notes in Computational Science and Engineering (LNCSE)}. {Springer-Verlag},
  Haidelberg, 2005.

\bibitem{saad_iterative_2003}
Y.~Saad.
\newblock {\em Iterative Methods for Sparse Linear Systems}.
\newblock Society for Industrial and Applied Mathematics, Philadelphia, {PA},
  {USA}, 2nd edition, 2003.

\bibitem{GPU:top500:2013}
{Top500}, 2013.
\newblock Available on line at: \url{http://www.top500.org/} (accessed on
  \today).

\end{thebibliography}
\bibliographystyle{abbrv}

\end{document}